\documentclass[12pt]{amsart}
\usepackage {amssymb,amsmath}
\makeatletter
\def\@cite#1#2{{\m@th\upshape\bfseries%
[{#1\if@tempswa{\m@th\upshape\mdseries, #2}\fi}]}} \makeatother
\theoremstyle{plain}
\newtheorem{thm}{Theorem}[section]
\newtheorem{lem}[thm]{Lemma}
\newtheorem{cor}[thm]{Corollary}

\theoremstyle{definition}
\newtheorem{rem}[thm]{Remark}

\newtheorem{defn}[thm]{Definition}

\newcommand{\Prf}{\noindent\textbf{Proof.\ }}
\newcommand{\bx}{\strut\hfill$\blacksquare$\medbreak}

\newcommand{\ca}{\mathrm{C}^*}

\newcommand{\bbC}{{\mathbb{C}}}

\newcommand{\bbN}{{\mathbb{N}}}

\newcommand{\bbT}{{\mathbb{T}}}
\newcommand{\bbZ}{{\mathbb{Z}}}
 \newcommand{\A}{{\mathcal{A}}}
 \newcommand{\B}{{\mathcal{B}}}

 \newcommand{\E}{{\mathcal{E}}}

\renewcommand{\H}{{\mathcal{H}}}

\renewcommand{\O}{{\mathcal{O}}}
\renewcommand{\P}{{\mathcal{P}}}
 
 \newcommand{\R}{{\mathcal{R}}}
\renewcommand{\S}{{\mathcal{S}}}

\newcommand{\upchi}{{\raise.35ex\hbox{$\chi$}}}
\newcommand{\qfor}{\quad\text{for}\quad}

\newcommand{\Alg}{\operatorname{Alg}}

\newcommand{\ran}{\operatorname{Ran}}

\newcommand{\spn}{\operatorname{span}}

\newcommand{\tr}{\operatorname{tr}}

\newcommand{\fix}{\operatorname{Fix}}
\newcommand{\sgn}{\operatorname{sgn}}
\newcommand{\Tab}{\operatorname{Tab}}
\newcommand{\Arr}{\operatorname{Arr}}
\newcommand{\SST}{\operatorname{SST}}
\def\bra#1{\langle #1|}
\def\ket#1{|#1 \rangle}
\def\one{{\mathchoice{\rm 1\mskip-4mu l}{\rm 1\mskip-4mu l}{\rm 1\mskip-4.5mu l}{\rm
1\mskip-5mu l}}}

\newcommand{\bofh}{\B(\H)}

\newcommand{\Ma}{{\mathbb M}}
\newcommand{\la}{\lambda}

\begin{document}

\baselineskip = 20pt plus 1pt minus 1pt

\title[Universal Collective Rotation Channels]%
{Universal Collective Rotation Channels and Quantum Error
Correction}
%
\author[M.Junge, P.T.Kim, D.W.Kribs]{Marius~Junge$^1$, Peter~T.~Kim$^2$, David~W.~Kribs$^{2,3,4}$}
\address{$^1$Department of Mathematics, University of Illinois at
Urbana-  \linebreak Champaign, Urbana, Illinois, USA 61801-2975}
\address{$^2$Department of Mathematics and Statistics, University of Guelph,
Guelph, ON, CANADA  N1G 2W1. }
\address{$^3$Institute for Quantum Computing, University of
Waterloo, Waterloo, ON, CANADA N2L 3G1. }
\address{$^4$Perimeter Institute for Theoretical Physics, 35 King
St. North, Waterloo, ON, CANADA N2J 2W9.}

\email{junge@math.uiuc.edu, pkim@uoguelph.ca,
\newline dkribs@uoguelph.ca}
%
\begin{abstract}
We present and investigate a new class of quantum channels, what
we call `universal collective rotation channels', that includes
the class of collective rotation channels as a special case. The
fixed point set and noise commutant coincide for a channel in this
class. Computing the precise structure of this $\ca$-algebra is a
core problem in a particular noiseless subsystem method of quantum
error correction. We prove that there is an abundance of noiseless
subsystems for every channel in this class and that the Young
tableaux combinatorial machine may be used to explicitly compute
these subsystems.
\end{abstract}
\maketitle

\section{Introduction}\label{S:intro}

The study of quantum channels is a central theme in quantum
computing and quantum information theory \cite{NC}.  A fundamental
class of quantum channels is known as the class of {\it collective
rotation channels}
\cite{BRS2,BRS,LPS,Fi,Lnoiseless,HKLP,HKL,KBLW,KLV,VKL,VFPKLC,ZL,Zan,Zar}.
This class has its roots in the postulates of quantum mechanics
and has recently played a key role in experimental efforts towards
realizing certain quantum error correction methods
\cite{Lnoiseless,VFPKLC}. Of particular interest in the current
study is the {\it noise commutant method of noiseless subsystems}.
This is a recently developed paradigm for passive quantum error
correction \cite{DG,Lnoiseless,HKL,IntroQEC,KLV,LCW,ZR}. In this
method, the `noise commutant' is used as a vehicle for encoding
states that are left immune to the noise of a given channel. The
operator algebras generated by such states are called `noiseless
subsystems'.

In this paper, we present a new class of quantum channels and
investigate them in the context of quantum error correction, with
specific reference to the noiseless subsystem method.  This class
is a generalization of the collective rotation class, which arises
as an important special case, hence we use the appellation
`universal collective rotation channels' to describe this class.
We prove that the noise commutant for every channel in this class
has rich structure and hence contains an abundance of noiseless
subsystems. To accomplish this, we use operator algebra techniques
to make an explicit connection with representation theory of the
symmetric group and, as a consequence, the Young tableaux
combinatorial machine \cite{FultonH,James} may be used to
explicitly compute these noiseless subsystems.

The paper is organized as follows. Section \ref{S:channels}
contains introductory material on quantum channels and quantum
error correction. In Section~\ref{S:universal} we  define and
establish basic properties of the class of universal collective
rotation ({\it ucr-}) channels.  We make the connection with
representation theory of the symmetric group in
Section~\ref{S:noisecomm} and show that the noise commutant for
ucr-channels is determined by a particular representation of the
symmetric group. In Sections~\ref{S:repnthy} and  \ref{S:examples}
we collect well-known facts from representation theory of the
symmetric group, with emphasis on Young tableaux combinatorics,
and work  through some low-dimensional examples. We finish with a
concluding remark in Section~\ref{S:conclusion} and discuss
possible avenues of further research.


\section{Quantum Channels and Noiseless Subsystems}\label{S:channels}


Let $\H$ be a (complex) Hilbert space and let $\bofh$ be the set
of bounded operators on $\H$.  When a basis for $\H$ is fixed and
$\dim\H  = k < \infty$, the algebra $\bofh$ may be identified with
the set of all complex $k\times k$ matrices $\Ma_k = \Ma_k
(\bbC)$. Throughout the paper, if we are given positive integers
$n\geq 1$ and $d\geq 2$, we let $\{ \ket{0}, \ket{1}, \ldots,
\ket{d-1} \}$ be a fixed orthonormal basis for $d-$dimensional
Hilbert space $\H_d = \bbC^d$ and  let $ \{ \ket{i_1 i_2 \cdots
i_n} : i_j\in\bbZ_d \}$ be the corresponding orthonormal basis for
$\H_{d^n} = (\bbC^d)^{\otimes n} $.

A linear map $\E : \bofh \rightarrow \bofh$ is {\it completely
positive} if for all $k \geq 1$ the ampliation maps $\one_k\otimes
\E : \Ma_k \otimes \bofh \rightarrow \Ma_k \otimes \bofh$  are
positive. See \cite{Kraustext,Paulsentext} for introductions to
the study of completely positive maps from different perspectives.
A {\it quantum channel} is a map $\E : \bofh \rightarrow \bofh$
that is completely positive and trace preserving. Given $\E$,
there is (\cite{Choi,Kraus}) a set of {\it noise operators}, or
{\it errors}, $\{ E_k \}$ on $\H$ such that
\begin{eqnarray}\label{channeldef}
\E(\rho) = \sum_{k} E_k \rho E_k^\dagger \qfor \rho \in\bofh.
\end{eqnarray}
Trace preservation means that the noise operators satisfy
\[
\sum_{k} E_k^\dagger E_k = \one,
\]
where $\one$ is the identity operator on $\H$. The channel is
unital if also,
\[
\E(\one) = \sum_{k} E_k E_k^\dagger  = \one.
\]



We will denote the fixed point set for $\E$ by
\[
\fix(\E) = \{ \rho \in \bofh : \E(\rho) = \rho \}.
\]
Further let $\A$ be the algebra generated by $\{E_k\}$ from
(\ref{channeldef}). This is called the {\it interaction algebra}
in quantum computing \cite{KLV}. It is a relic of the channel in
the sense that the same algebra is obtained whatever the choice of
noise operators in (\ref{channeldef}). This is most succinctly
seen in the case of a unital channel. In general, $\fix(\E)$ is
just a $\dagger$-closed subspace of $\bofh$, but in the case of a
unital channel $\E$, the so-called {\it noise commutant}
\[
\A^\prime = \{ \rho \in \bofh : \rho E_k =E_k \rho,  \,\,\forall k
\}
\]
coincides with the fixed point  set \cite{BS,Kchannel}:
\[
\fix(\E) = \A^\prime.
\]
In particular, $\fix(\E) = \A^\prime$ is a $\dagger$-closed
operator algebra (a finite dimensional $\ca$-algebra
\cite{Arvinvite,byeg,Tak}). In this case the von Neumann double
commutant theorem from operator algebras shows how the algebra $\A
= \A^{\prime\prime} = \fix(\E)^\prime$ only depends on the
channel.

Every finite dimensional $\ca$-algebra is unitarily equivalent to
an orthogonal direct sum of `ampliated' full matrix algebras; {\it
i.e.}, there is a unitary operator $U$ such that
\[
U \A U^\dagger = 
\bigoplus_{k=1}^r \,\,\big( \one_{m_k}\otimes \Ma_{n_k}  \big).
\]
From the representation theory perspective, a factor $\one_{m_k}
\otimes \Ma_{n_k}$ corresponds to an $n_k$-dimensional irreducible
representation appearing with multiplicity $m_k$. With this form
for $\A$ given, the structure of the commutant up to unitary
equivalence is easily computed by
\begin{eqnarray}\label{commform}
U \fix(\E) U^\dagger = U\A^\prime U^\dagger = 
\bigoplus_{k=1}^r\,\, \big( \Ma_{m_k} \otimes \one_{n_k} \big).
\end{eqnarray}
(See \cite{Lnoiseless,HKL,KBLW,VKL,ZL,Zan} for  more detailed
discussions in connection with quantum information theory.)

Given a quantum channel $\E$ with noise operators $\{E_k\}$, the
structure of the noise commutant $\A^\prime$  can be used to
prepare density operators for use in the noiseless subsystem
method of error correction. This is a passive method of quantum
error correction, in the sense that such operators will remain
immune to the effects of the noise  of the channel, without active
intervention.
Thus, computing the precise structure of  $\A^\prime$ as in
(\ref{commform}) is of fundamental importance in this method. We
mention that for experimental reasons \cite{Rayprivate},  only one
matrix algebra $\Ma_{m_k} \otimes \one_{n_k}$ may be used at a
time in this manner. Hence it is also desirable to find the
largest full matrix algebra which is a subalgebra of the noise
commutant.

\section{Universal Collective Rotation Channels}\label{S:universal}

For the rest of the paper, given a positive integer $d\geq 2$ we
write $\Ma_d$ for the operator algebra $\B(\bbC^d)$ represented as
$d\times d$ complex matrices with respect to the standard basis
$\{ \ket0 , \ldots, \ket{d-1}\}$ for $\bbC^d$. Further let
$\Ma_{d,sa}$ be the subset of self-adjoint matrices inside
$\Ma_d$.

Fix $n\geq 1$. Given $1\leq k \leq n$ we define a representation
of $\Ma_d$ on $\H_{d^n}$ by
\[
\omega_{k}(x) = \one_d \otimes \cdots \otimes \one_d \otimes
 \underbrace{x}_{\mbox{\scriptsize $k$-th position}}\otimes \one_d
 \otimes
 \cdots \otimes \one_d
\]
for all $x\in  \Ma_{d}$. Then we may define sums of independent
copies of $x$ by
\[
 u_n(x) =  \sum_{k=1}^n \omega_k(x) \qfor x\in\Ma_d .
\]

\begin{defn}
Given a finite subset $\S \subset \Ma_{d,sa}$, we define a {\it
universal collective rotation (ucr-) channel} $\E_\S$ by
\[
\E_\S (\rho) = \frac{1}{\sqrt{|\S|}} \sum_{x\in\S} e^{i\theta_x
u_n(x)}\,\, \rho \,\, e^{-i\theta_x u_n(x)} \qfor \rho \in
\B(\H_{d^n}),
\]
where $\{ \theta_x :x\in\S\} $ are non-zero angles.
\end{defn}

Given a set of operators $\R$, define $\Alg \R$ to be the operator
algebra generated by $\R$. This is the set of all polynomials in
the elements of $\R$. When $\R$ is a self-adjoint set, $\Alg \R$
is a $\ca$-algebra. Through a standard functional calculus
argument from operator theory, it follows that the interaction
algebra $\A_\S$ for $\E_\S$ is obtained as $\A_\S \equiv \Alg \{
e^{i\theta_x u_n(x)} : x\in\S\} = \Alg \{ u_n(x) : x\in\S\}$. Thus
by von Neumann's double-commutant identity we have
\[
\A_\S = \{ e^{i\theta_x u_n(x)} : x\in\S\}^{\prime\prime} = \{
u_n(x) : x\in\S\}^{\prime\prime}.
\]
Notice that $\A_\S$ is independent of the choice of (non-zero)
angles $\theta_x$. As an application of the fixed point theorem
from \cite{BS,Kchannel} we obtain the following.

\begin{thm}
If $\S$ is a finite subset of $\Ma_{d,sa}$, then the ucr-channel
$\E_\S$ satisfies
\[
\fix(\E_\S) = \A_\S^\prime.
\]
\end{thm}

Observe that $\one_d$ belongs to  $\A_\S$ from its
characterization as a  bicommutant. Since
$u_n(\one_d)=n\one_{d^n}$, we may always add $\one_d$ to $\S$
without changing the properties of $\fix(\E_S)$. This motivates
the following definition.

\begin{defn}
We  will say that $\S$ is {\it maximal} if $\spn \{x:x\in \S\}$
contains all matrices $x\in\Ma_d$ with $\tr(x)=0$.
\end{defn}

It turns out that for maximal $\S$ the algebra $\A_\S$ is a
well-known object in representation theory.

\begin{rem}
To place the class of ucr-channels in context, we note that the
ucr-channels for $d=2$ and general $n$ are the class of
`two-level' collective rotation channels from quantum computing
\cite{BRS2,BRS,LPS,Fi,Lnoiseless,HKLP,HKL,KBLW,KLV,VKL,VFPKLC,ZL,Zan,Zar}.
The noise operators in this case are also denoted by $J_x,J_y,J_z$
and they arise in quantum mechanics as the  canonical
representation of the angular momentum relations \cite{QM12}. From
the noiseless subsystem/quantum error correction perspective, the
algebra $\fix(\E_\S) = \A_\S^\prime$ for this subclass of
ucr-channels, and natural $d$-dimensional representations of the
$J_k$ operators, has been analyzed in \cite{HKLP} from an operator
theory cum quantum mechanics point of view.
\end{rem}

\section{Representation Theory and the Noise Commutant}\label{S:noisecomm}

In this section we identify the structure of the noise commutant
in terms of representation theory for the symmetric group. We
begin with some notation. We shall denote the $n$-fold tensor
product of $\Ma_d$ by
\[
\Ma_d^{\otimes n} =  \underbrace{\Ma_d\otimes \cdots \otimes
\Ma_d}_{\mbox{\scriptsize $n$-times}} \cong \Ma_{d^n}.
\]
Let ${\rm Sym}^n \Ma_d$ be the subalgebra of $\Ma_{d^n}$ generated
by the symmetric tensor products; that is, ${\rm Sym}^n \Ma_d$ is
the algebra generated by the operators
\[
\Phi_n (x_1\otimes \cdots \otimes x_n) = \frac{1}{n!} \sum_{\pi\in
S_n} x_{\pi(1)} \otimes \cdots \otimes x_{\pi(n)},
\]
where each $x_i \in \Ma_d$ and $S_n$ is the permutation group on
$n$ letters.

In terms of representation theory, we may equally well consider
the representation $\pi: {\rm GL}(d) \rightarrow {\rm GL}(d^n)$ given by
$\pi(u) = u\otimes \cdots \otimes u$, and then we have
\[
{\rm Sym}^n \Ma_d = \pi({\rm GL}(d))'',
\]
where ${\rm GL}(d)$ is the group of $d\times d$ nonsingular complex matrices.
This tensor product representation of ${\rm GL}(d)$ is in `duality' with
the representation of the symmetric group $S_n$ defined on vector
tensors by
\[
\pi(\sigma) (h_1\otimes \cdots \otimes h_n) = h_{\sigma(1)}
\otimes \cdots \otimes h_{\sigma(n)},
\]
for $\sigma\in S_n$ and $h_1, \ldots ,h_n\in\H_d$. In this
context, Schur's classical duality theorem  reads as

\begin{thm}\label{schur}
$\pi(S_n)' = {\rm Sym}^n \Ma_d$.
\end{thm}

We use the following characterization of ${\rm Sym}^n \Ma_d$
below.

\begin{lem}\label{symmetric}
For positive integers $d$ and $n$, we have
\[
{\rm Sym}^n \Ma_d =  \{ x^{\otimes n} : x\in\Ma_d \}''  = \{ u_n
(x) : x\in\Ma_d \}''
\]
\end{lem}

\Prf It is clear that ${\rm Sym}^n \Ma_d$ contains the
$\ca$-algebra
\begin{eqnarray*}
\B =  \{x^{\otimes n} : x\in\Ma_d\}''  &=& \spn \{ x^{\otimes n} :
x\in\Ma_d\} \\ &=& \Alg \{ x^{\otimes n} : x\in\Ma_d\}
\end{eqnarray*}
as a subalgebra. For the converse inclusion, let $x_1, \ldots
,x_n \in \Ma_d$ and consider the complex matrix integral
\begin{align*}
& \int_{z_1,\ldots ,z_n\in\bbT} \Big( \sum_{j=1}^n z_j
x_j\Big)^{\otimes n} \,\,\frac{dz_1}{z_1} \cdots \frac{dz_n}{z_n} \\
& = \sum_{j_1,\ldots ,j_n =1}^n \left( \int_{z_1,\ldots
,z_n\in\bbT} \prod_{r=1}^n z_{j_r} \frac{dz_1}{z_1} \cdots
\frac{dz_n}{z_n} \right) (x_{j_1} \otimes \cdots \otimes x_{j_n})
\\
& = \sum_{j_1,\ldots ,j_n =1}^n \left( \prod_{s=1}^n
\int_{z\in\bbT} z^{|\{r:j_r = s\}|} \frac{dz}{z} \right) (x_{j_1}
\otimes \cdots \otimes x_{j_n})
\\
& =(2\pi i)^n \sum_{\pi\in S_n} x_{\pi(1)} \otimes \cdots \otimes
x_{\pi(n)},
\end{align*}
where $\bbT$ denotes the unit circle in the complex plane. It
follows that $\Phi_n (x_1\otimes \cdots \otimes x_n)$ belongs to
$\B$ for any choice of $x_1, \ldots, x_n$, and hence $\B$
coincides with ${\rm Sym}^n \Ma_d$.

On the other hand, it is clear by definition that ${\rm Sym}^n
\Ma_d$ contains the algebra $\{ u_n(x):x\in\Ma_d\}''$ generated by
the $u_n(x)$. Moreover, a consideration of the expansion for
$u_n(x)^n$ shows that $x^{\otimes n}$ belongs to this double
commutant for all $x\in\Ma_d$. For the sake of brevity let us
observe this fact for $n=2$ and $n=3$:
\begin{eqnarray*}
x\otimes x & =& \frac{1}{2!} \Big( u_2(x)^2 - u_2(x^2)\Big) \\
x\otimes x \otimes x & =& \frac{1}{3!} \Big( u_3(x)^3 - 3
u_3(x^2)u_3(x) - 2 u_3(x^3) \Big).
\end{eqnarray*}
In fact, for all $x\in\Ma_d$, the tensor product $x^{\otimes n}$
belongs to the algebra $\Alg \{ u_n(x^p) : 1\leq p \leq n\}$. Thus
the second characterization of ${\rm Sym}^n \Ma_d$ follows. \bx

Observe that as a consequence of this proof, we also have ${\rm
Sym}^n \Ma_d =  \{ u_n(x) : x\in\Ma_{d,sa} \}''.$  We can now
explicitly link the noise commutant for these channels with
representation theory of the symmetric group.

\begin{thm}\label{mainthm}
Let $\S\subset \Ma_{d,sa}$ be a maximal system, then
\[
\fix (\E_\S) = \A_\S^\prime = \pi (S_n)''.
\]
Moreover, for an arbitrary finite set $\S\subset \Ma_{d,sa}$, we
have
\[
\fix (\E_\S) \supseteq \pi (S_n)''.
\]
\end{thm}

\Prf If $\S$ is maximal, then the interaction algebra
$\fix(\E_\S)' = \A_\S = \{ u_n(x) : x\in\S\}'' = \{ u_n(x) :
x\in\Ma_d\}''$ coincides with $\pi(S_n)'$ by
Lemma~\ref{symmetric}. For the second assertion, a given finite
subset $\S \subset \Ma_{d,sa}$ is contained inside a maximal
system $\S_{\max}$. Hence $\A_\S \subseteq \A_{\S_{\max}}$ and
\[
\fix(\E_\S) = \A_\S^\prime \supseteq \A_{\S_{\max}}^\prime = \pi
(S_n)''.
\]
\bx

\section{Computing Noiseless Subsystems Via Young Tableaux}\label{S:repnthy}

In this section, we collect well-known facts from the
representation theory of $S_n$ that allow us to describe $\fix
(\E_\S) = \pi (S_n)''$ in an explicit manner. Recall that this is
imperative for using the structure of the noise commutant to
produce noiseless subsystems.

For the discussion in this section, we shall fix positive integers
$d\geq 2$ and $n\geq 2$. Let $\{\ket0, \ldots, \ket{d-1}\}$ be the
orthonormal basis for $\H_d$ corresponding to a given $d$-level
quantum system, and let
\[
\big\{ \ket{i_1\cdots i_n} : 0\leq i_j < d, \, 1\leq j \leq
n\big\}
\]
be the corresponding basis for $\H_{d^n}$. Observe that the set of
$n$-tuples $\{i_1,\ldots, i_n\}$ is in one-to-one correspondence
with the set of functions $f:\{1,\ldots,n\}\to \{0,\ldots,d-1\}$.
So we may define functions $k_l$ for $0 \leq l < d$ by
\[
k_l(i_1,....,i_n) =  \#\{1\leq j\leq n \,\,| \,\,i_j=l\} \qfor
0\leq i_j < d,
\]
and we have   $\sum_{l=0}^{d-1} k_l(i_1,...,i_n)=n$.

Now, given positive integers $k_0, \ldots ,k_{d-1}$ with each
$0\leq k_l \leq n$, we define a corresponding subspace of
$\H_{d^n}$ by
\[
\H_{k_0,....,k_{d-1}} =  \spn\big\{ \ket{i_1\cdots i_n} :
 k_l(i_1,....,i_n)=  k_l,\, 0\leq l < d \big\} .
\]
Notice that $\H_{d^n} = \bigoplus \H_{k_0,\ldots,k_{d-1}}$, where
the direct sum runs over all choices of $k_0,\ldots,k_{d-1}$. Clearly,
$\H_{k_0,....,k_{d-1}}$ is an invariant (hence reducing) subspace
for the action of the symmetric group $S_n$. More importantly, the
irreducible subspaces, or equivalently the decomposition factors
of $\H_{k_0,...,k_{d-1}}$ are completely characterized. The key
ingredient in this characterization is the notion of Young
tableaux.

Given $\lambda_1 \geq \la_2 \geq \ldots \geq \lambda_r$, a
non-increasing sequence of positive integers  with
$\sum_i\la_i=n$, put $\la=(\la_1,\ldots ,\la_r)$. Then the
associated {\it $\la$-diagram} is defined as
 \[
 [\la] = \{c_{ij} : 1\leq i\leq r, \, 1\leq j\leq \la_i\},
 \]
where $c_{ij}$ denotes a `cell' in $[\la]$. Simply put, $[\la]$ is
a diagram with $d$ rows of cells which are left justified and
$\la_i$ cells in the $i$th row. A {\it $\la$-tableau} is a
bijective function $t:[\la]\to \{1,\ldots,n\}$. Clearly, $S_n$
acts by composition $\sigma t = \sigma\circ t$ on $\la$-tableaux.
Given a $\la$-tableau, the {\it column stabilizer} $C_{t}$ is the
subgroup of $S_n$ which leaves the columns of $\la$ setwise fixed.
Similarly, the {\it row stabilizer} $R_{t}$ is the subgroup of
$S_n$ which leaves the rows of $\la$ setwise fixed. Two tableaux
$t_1$ and $t_2$ are {\it equivalent} if there exists a permutation
$\sigma\in R_{t_1}$ such that $\sigma t_1=t_2$. In particular,
this means that the set of {\it tabloids} $\Tab_{\la}=\{\{t\}: t\
\mbox{a\  $\la$-tableau}\}$ of equivalence classes is indexed by
all partitions $(A_1,...,A_r)$ of $\{1,\ldots,n\}$  such that the
cardinalities $|A_1|=\la_1$,\ldots,$|A_r|=\la_r$.   Given a
$\la$-diagram $[\la]$, consider the $(i,j)$-cell $c_{ij}$ in
$[\la]$. The {\it hook length} $h(i,j)$ for $c_{ij}$ is the number
of cells directly below $c_{ij}$ in the $j$th column of $[\la]$
plus the number of cells to the right of $c_{ij}$ in the $i$th row
of $[\la]$ plus one (for the cell $c_{ij}$ itself). Formally,
 \[
 h(i,j)= \la_i+\la_j'+1-i-j  ,
 \]
where $\la_j'$ is the number of elements in the $j$-th column.
Also recall that a tableau $t:[\la]\to \{1,\ldots,n\}$ is {\it
standard} if the numbers increase along rows and increase down
columns. The abstract $S_n$-module that has an orthonormal basis
in bijective correspondence with elements of $\Tab_{\la}$ is
denoted by $M^\lambda$, so that
\[
M^\lambda = \spn \{ e_{\{t\}} : \{ t \} \in \Tab_{\la} \}.
\]
The {\it Specht module} $\S^{\la}$ is the submodule of $M^{\la}$
generated by the `polytabloids'
\[
e_t\equiv \sum_{\sigma\in C_t}(\sgn \sigma) e_{\sigma \{t\}} \in
M^{\la}.
\]
Let us summarize the following facts (see chapter 7 in
\cite{James}).

\begin{thm}\label{spechtthm} Let $k_0,\ldots,k_{d-1}$ be positive integers and
consider a partition of $\{1,\ldots,n\}$ into sets
$A_0,\ldots,A_{d-1}$ with $|A_l|=k_l$. Let $\la$ be the
non-increasing rearrangement of $(k_0, \ldots, k_{d-1})$.  Then
$\H_{k_0,\ldots,k_{d-1}}$ is isomorphic as an $S_n$-module to
$M^{\la}$.

Every polytabloid $e_t$ is a cyclic vector for the irreducible
module $\S^{\la}$.  The dimension of $\S^{\la}$ is given by the
`hook length formula'
 \[
 \dim \S^{\la} = \frac{n!}{\prod \mbox{hook lengths in
 $[\la]$}} ,
 \]
and a basis for $\S^{\la}$ is given by
 \[
 \{e_t: t \mbox{ standard $\la$-tableau } \} .
 \]

Finally, every finite dimensional irreducible representation of
$S_n$ is unitarily equivalent to a Specht module representation
$\pi_{\la}$, where $\pi_{\la}$ is the representation of $S_n$ on
$\S^\la$ defined by $\pi_\la (\sigma) e_t \equiv e_{\sigma t}$.
\end{thm}

For the next discussion let us fix numbers $(k_0,\ldots,k_{d-1})$
and let us denote by $\mu=(\mu_0,\ldots,\mu_{d-1})$ the
non-increasing rearrangement of $(k_0,\ldots,k_{d-1})$.  The
$S_n$-module $\H_{k_0,...,k_{d-1}}$ decomposes into a direct sum
of irreducible submodules. Fortunately, these submodules and their
multiplicity are completely characterized by Young's rule.
Moreover, below we shall describe how the decomposition into
irreducible submodules of $\H_{k_0,...,k_{d-1}}$ is related to,
and determined by, the decomposition of $\H_\mu \equiv
\H_{\mu_0,...,\mu_{d-1}}$. (This allows us to explicitly identify
links between irreducible subspaces for the representation $\pi$.)
Here the key combinatorial tool is the notion of a semistandard
tableau.

We generalize the notion of $\la$-tableau, by saying that
$T:[\la]\to \bbN$ is a {\it $\la$-tableau of type
$\mu=(\mu_0,....,\mu_{d-1})$} if
 \[
 \#\big\{c_{ij}: T(c_{ij})=l\big\} = \mu_l \qfor l=0,\ldots,d-1.
 \]
Then $T$ is called {\it semistandard} if the numbers that $T$
assigns to the cells of the diagram determined by $\la$ are
non-decreasing along rows and strictly increasing down columns.
Let us fix a bijection $t_0:[\la]\to \{1,\ldots,n\}$. Then $S_n$
acts on the sets $I(\la,\mu)$, the set of $\la$-tableau of type
$\mu$, {\it via}
 \[
 \sigma(T) =  T \, t_0^{-1}\sigma \, t_0  \qfor \sigma\in S_n .
 \]

Given $t_0$, we will say that $T_1$ and $T_2$ are {\it row} ({\it
column}) equivalent, and write $T_1 \sim_{t_0}^r T_2$, if $\sigma
T_1=\sigma T_2$ holds for all permutations $\sigma$ in the row
(respectively column) stabilizer of $t_0$. In particular, this
means that $T_1$ and $T_2$ are row equivalent if and  only if
$T_2$ is obtained from $T_1$ by permuting the entries in each row
accordingly.

In order to define the linking module maps we first need an
appropriate bijection. We denote by $\P_{\mu_0,\ldots,\mu_{d-1}}$
the set of partitions $(A_0,...,A_{d-1})$ of $\{1,\ldots ,n\}$
such that $|A_l|=\mu_l$. Then  $\P_{\mu_0,\ldots,\mu_{d-1}}$
induces a natural relabelling of the standard basis for $\H_{\mu}$
by
 \begin{eqnarray}\label{identification}
 f_{A_0,...,A_{d-1}} =  \ket{i_1 \cdots i_n}
 \end{eqnarray}
where $A_l=\{1\leq j\leq n \,\,|\,\, i_j=l\}$ for $0\leq l < d$.
(Every $n$-tuple $(i_1,\ldots,i_n)$ is associated with a unique
$d$-tuple of sets $(A_0,\ldots,A_{d-1})$ defined in this way.)

Next we define $ \gamma_{t_0}:I(\la,\mu)\to
\P_{\mu_0,...,\mu_{d-1}}$ by
\[
\gamma_{t_0}(T)= (A_0, \ldots , A_{d-1}),
 \]
where
\[
A_l = \{ 1\leq j \leq n \,\, | \,\, T t_0^{-1}(j) = l \} \qfor
0\leq l < d.
\]
Every $\la$-tableau  $T$ of type $\mu$ induces an $S_n$-module map
$\Theta_{T}:M^{\la}\to M^{\mu}$ by
 \[
 \Theta_T(e_{\{t_0\}})= \sum_{T'\sim_{t_0}^r T, \,\,\,\gamma_{t_0}(T')=(A_0,\ldots,A_{d-1})}
    f_{A_0,...,A_{d-1}}  .
 \]
Clearly this extends to an $S_n$-module homomorphism by defining
 \[
 \Theta_T(e_{\{\sigma(t_0)\}}) = \sigma(\Theta_T(e_{\{t_0\}}))  .
 \]

This rather abstract description is in fact very concrete. Given
indices $i_1,\ldots,i_n\in \{0,\ldots , d-1\}$ and a $\la$-tableau
$t:[\la]\to \{1,\ldots,n\}$ we form the generalized tableau
$t_{\ket{i_1 \cdots i_n}}:[\la]\to \{0,\ldots ,d-1\}$ by
 \[
 t_{\ket{i_1 \cdots i_n}}(c_{ij})= i_{t_0(c_{ij})} .
 \]
This means we write the entries $i_1,\ldots,i_n$ into $\la$
following the order given by $t_0$. Then we say that
 \[
 (i_1,\ldots,i_n)\sim_{t_0} (i_1',\ldots,i_n')
 \]
if there exists a permutation $\sigma \in S_n$  such that
 $i'_j= i_{\sigma(j)}$ for $1\leq j \leq n$ and
 $t_0^{-1}\sigma t_0$ leaves  the rows of $\la$ invariant.
Therefore, we obtain
 \[
 \Theta_T(e_{\{t_0\}})= \sum_{(i_1,\ldots,i_n)\sim_{t_0} \gamma_{t_0}(T)} \ket{i_1\cdots i_n} ,
 \]
where here we identify $\gamma_{t_0}(T)$ with the $n$-tuple
determined by the partition $\gamma_{t_0}(T) = (A_0, \ldots
,A_{d-1})$ as in (\ref{identification}).

For example, let $d=3$, $n=5$ and let $t_0:[\la]\to
\{1,\ldots,5\}$ be given by
 \[
 t_0 = \begin{array}{|c|c|c|c|} \hline 1 & 3 & 2& 4  \\ \hline
                                 5 \\ \cline{1-1}
                                 \end{array}
 \]
and $T:[\la]\to \{0,1,2\}$ be given by
 \[
 T = \begin{array}{|c|c|c|c|} \hline 0 & 0 & 1& 1 \\ \hline
                                 2 \\ \cline{1-1}
                                 \end{array} .
 \]
This yields, by reading off the entries from the corresponding
position in the diagram,
 \[
\gamma_{t_0}(T)= (A_0, A_1,A_2) = (\{1,3\}, \{2,4\}, \{5\}).
 \]
Following (\ref{identification}), $\gamma_{t_0} (T)$ is identified
with $(i_1,i_2,i_3,i_4,i_5) = (0,1,0,1,2)$. Moreover, the list of
equivalent indices is:
  \begin{eqnarray*}
  \Big\{ (0,1,0,1,2), (0,1,1,0,2),(0,0,1,1,2), (1,0,0,1,2), \\
  (1,0,1,0,2),(1,1,0,0,2)\Big\} .
  \end{eqnarray*}
Indeed, according to $t_0$ we have to fix the $5$th coordinate and
the other four vary in all possible ways.  Thus we have
 \[
 \Theta_T(e_{\{t\}})= \sum_{(i_1,\ldots,i_n)\sim_{t}\gamma_t(T)}
 \ket{i_1\cdots i_n} \qfor e_{\{t\}}\in M^\la .
 \]

Following Young's rule (see chapter 2, \cite{James}) we obtain:

\begin{thm}
Let $\mu=(k_0^*,\ldots,k_{d-1}^*)$ be the non-increasing
rearrangement of $(k_0,\ldots,k_{d-1})$. Let
$\la=(\la_1,\ldots,\la_r)$ be such that $\lambda_1 \geq \ldots
\geq \la_r$ and $\sum_i \la_i=n$. Then
 \[
 \H_{k_0,\ldots,k_{d-1}}^{\la} \equiv \spn\left\{ \Theta_{T}\Big(\sum_{\sigma \in C_t} (\sgn \sigma)e_{\sigma\{t\}} \Big) :
   T\in I(\la,\mu) , \,\,
 t  \mbox{ $\la$-tableau }\right\}
 \]
is an irreducible $S_n$-submodule. The restriction of the
representation $\pi$ to $\H_{k_0,\ldots,k_{d-1}}^{\la}$ is
equivalent to the irreducible representation $\pi_\mu$ of $S_n$ on
$S^{\mu}$ and has multiplicity
  \[
  m = \#\big\{ T : T \mbox{ semistandard $\la$-tableau of type } \mu \big\}  .
  \]
\end{thm}

If we collect all this information for  all $(k_0,...,k_{d-1})$,
we can describe the full representation $\pi$ of $\bbC[S_n]$:

\begin{cor}
Let $\la=(\la_1,\ldots,\la_r)$ be such that $\lambda_1 \geq \ldots
\geq \la_r$ and $\sum_{i}\la_i=n$, and let $P_{\la}$ be the
projection of $\H_{d^n}$ onto
 \[
 \H^{\la}\equiv   \bigoplus_{k_0,\ldots ,k_{d-1}}\,\, \H_{k_0,\ldots,k_{d-1}}^{\la}   ,
 \]
where the sum indexes over all $k_0,\ldots ,k_{d-1}$ such that
$\sum_{l=0}^{d-1} k_l =n$.

Then $P_\la$ is the minimal central projection for $\pi(S_n)''$
which supports the irreducible submodule $\S^\la$. Moreover,
$P_{\la}\pi P_{\la}$ is equivalent to the representation
$\pi_{\la}$ on $\S^{\la}$ with multiplicity
 \[
 m_{\la,d} = \sum_{\mu_0\geq \cdots \geq \mu_{d-1}} \Arr(\mu)
 \SST(\mu),
 \]
where
 \[
 \Arr(\mu)=\#\Big\{(k_0,...,k_{d-1}):
 (k_0^*,...,k_{d-1}^*)=(\mu_0,...,\mu_{d-1})\Big\}
 \]
and
\begin{eqnarray*}
\SST(\mu) = \#\Big\{ T:[\la]\to \{0,\ldots,d-1\} &\Big|& T \mbox{
semistandard} \\ & & \mbox{ $\la$-tableau of type }
 \mu \Big\}.
\end{eqnarray*}

In particular, for a maximal system $\S$,
 \[
 \fix(\E_S) = \A_\S^\prime = \pi(S_n)'' \cong  \sum_{m_{\la,d}\neq 0}
 \Ma_{\dim(\S^{\la})}\otimes \one_{m_{\la,d}}
 \]
describes the representation in irreducible parts with
multiplicity.
\end{cor}

\begin{cor}
Let $\S$ be a maximal system, then $\A_\S$ is
isomorphic to
 \[
 \A_\S\cong \sum_{m_{\la,d}\neq 0} \Ma_{m_{\la,d}}\otimes
 \one_{\dim \S^\la}
 \]
and the multiplicity of the component $\Ma_{m_{\la,d}}$ is given
by $\dim \S^{\la}$.
\end{cor}

Let us mention that $\H^{\la}$ may also be described by the
so-called Garnier relations. Given $\la=(\la_1,\ldots,\la_r)$, we
fix the tableau $T_{\la}$ such that $T_{\la}(c_{kj})=k$ for all
cells $c_{kj}$ in $[\la]$. It follows that every index
$i=(i_1,\ldots,i_n)$ defines a tableau $T_i:[\la]\to
\{0,\ldots,d-1\}$ given by $T_i(c_{kj})= i_{T_{\la}(c_{kj})}$. Let
$G(J)$ be the collection of coset representatives $\{\nu X:\nu \in
Y\}$, where $Y$ is the subgroup of $S_n$ which fixes every element
outside both $C_h(T_{\la})\cup J$ and $Y=X\cap C(T)$. Then as is
proved in  \cite{Green}(p.66, 5.2b), $\ket{\psi}\in \H^{\la}$ if
and only if
\begin{enumerate}
 \item $\bra{\psi}\ket{i}=0$ for all $i$ such that $T_i$ has equal entries
 in two distinct places in the same column.
 \item $\pi(\sigma)(\ket{\psi})= \sgn(\sigma)\ket{\psi}$ for all $\sigma$ in the column
 stabilizer of $T_{\la}$.
 \item $\sum_{\nu\in G(J)} \sgn(\nu)\pi(\nu^{-1})\ket{\psi}=0$ for any
 non-empty set  in the column stabilizer of $C_{h+1}(T_{\la})$.
\end{enumerate}

\section{Examples}\label{S:examples}

\subsection{The case $d=2$ and general $n$}

As mentioned above, the case of $d=2$ and general $n$  was
extensively examined in \cite{HKLP}. Let us  indicate how this can
be accomplished with Young tableaux.

When $d=2$, we have the pairs $(k_1,k_2)$ given by $(n-k,k)$ where
$k=0,...,n$. In terms of $\la$-tableau we have to calculate
$m_{\la,2}$. In terms of types we only have to consider diagrams
$\mu_k=(n-k,k)$ where $0\leq k\leq 2n$.  But we have to be aware
that every type allows  combinations $(n-k,k)$ and $(k,n-k)$.
Given $\la=(\la_1,...,\la_r)$, we observe that to obtain a
semistandard tableau, we must have $r=2$. Indeed, we are forced to
put  $0$'s in the first row on the first $\la_2$ positions and
$1$'s in the second row. Thus for fixed $k,j$ with $2k\leq n$ an
$2j\leq n$, we need $k\leq j$ in order to produce a $\la$-tableau
of type $\mu$. Since, we also know that there are $n-k$ $0$'s, we
do not have a choice and we have to put them all in the first row
one after another. Thus for a fixed $\la$, we find
 \[
 m_{(n-j,j),2} = \sum_{k\leq j}^{\lfloor n/2 \rfloor} 2 + 1
 = 2  (\lfloor n/2 \rfloor -j)+1
 \]
if $n$ is even and
 \[
 m_{(n-j,j),2} = 2 (\lfloor n/2 \rfloor -j)
 \]
if $n$ is odd, where $\lfloor \cdot \rfloor$ denotes the greatest
integer part of some number.

We also have to understand $\dim(\S^{\la})$. If $\la=0$, we get
$\dim(\S^{(n,0)})=1$. If $1\leq j<\frac{n}{2}$, we see for cells
$c_{1l}$ with $l\leq j$ the hook length is $1+(n-j+1-i)$. This
yields $n (n-2j+1)/j! (n-j+1)!$ and hence
 \[
 \dim(\S^{(n-j,j)})= \begin{cases} 1 & \mbox{ if } j=0  \\
  \frac{n-2j+1}{n+1} {n+1 \choose j} & \mbox{ if }
  1<j\leq \frac{n}{2}
  \end{cases}
   .\]

Let us consider the examples $n=4$ and $n=5$. Then
 \[ \dim \S^{(4,0)}= 1, \quad \dim \S^{(3,1)} = 3 , \quad
 \dim \S^{(2,2)} =  2 \]
and
 \[ m_{(4,0),2}= 5, \quad m_{(3,1),2} = 3  ,  \quad m_{(2,2),2} = 1  .\]
In the case $n=5$, we have
 \[ \dim \S^{(5,0)}= 1, \quad \dim \S^{(4,1)}= 4  , \quad
 \dim \S^{(3,2)} =  5 \]
and
 \[ m_{(5,0),2}= 6, \quad m_{(4,1),2} = 4 , \quad m_{(3,2),2} = 2  .\]
Bases for $\H_{2^n}$ which yield the associated algebra
decompositions may be computed as well. Below we do this for a
more intricate example.

\subsection{The case $d=3$ and $n=4$}

If $d=3$ and $n=4$, then the set of $\la$-diagrams which admit
semistandard tableaux is given by
 \begin{eqnarray*}  \left\{(4)=\begin{array}{|c|c|c|c|} \hline
 x&x&x&x \\ \hline \end{array} \quad
 (3,1)=\begin{array}{|c|c|c|} \hline x&x&x\\ \hline
                          x \\ \cline{1-1}  \end{array} \quad
  (2,2)=\begin{array}{|c|c|} \hline x&x\\ \hline
                          x&x \\ \cline{1-2} \end{array} \right. \\
\left.  (2,1,1)=\begin{array}{|c|c|} \hline x&x\\ \hline
                            x\\ \cline{1-1}
                            x \\ \cline{1-1} \end{array} \right\}
                            .\end{eqnarray*}
In this case, $\pi(S_4)$ acts on $\H_{d^n} = \H_{81}$.
As in Theorem~\ref{spechtthm}, $M^{(4)}$ is isomorphic to
$\H_{4,0,0}$, $\H_{0,4,0}$ and $\H_{0,0,4}$; $M^{(2,2)}$ is
isomorphic to $\H_{2,2,0}$, $\H_{2,0,2}$ and $\H_{0,2,2}$;
$M^{(2,1,1)}$ is isomorphic to $\H_{2,1,1}$, $\H_{1,2,1}$ and
$\H_{1,1,2}$; {\it etc}, so that the multiplicities for the
$M^\lambda$ are $3$ for $M^{(4)}$, $M^{(2,2)}$ and $M^{(2,1,1)}$
and $6$ for $M^{(3,1)}$. The dimensions of the Specht modules
$\S^{\la}$ using the hook length formula are given by
 \[ \dim \S^{(4)}= 1,\,\,\,\, \dim \S^{(3,1)} = 3 ,\,\,\,\,
 \dim \S^{(2,2)} =  2   ,\,\,\,\, \dim \S^{(2,1,1)}= 3 .\]

Now, we have to compute the multiplicities  of $\S^{\la}$ in
$M^{\mu}$. If $\mu=(4) (=(4,0,0))$, then the only semistandard
tableau of type $\mu$ is $\la = (4)$ with $0$ in each cell. Thus
$M^{(4)}\cong \S^{(4)}$. Further, every $M^\mu$ supports a single
copy of $\S^{(4)}$ {\it via} the $\la =(4)$-tableau with cell
entries given by $\mu$. For $\mu=(3,1)$, the possible semistandard
tableaux are
 \[ \begin{array}{|c|c|c|c|} \hline 0&0&0&1 \\ \hline \end{array}  \quad\quad
 \begin{array}{|c|c|c|} \hline 0&0&0\\ \hline
                    1 \\ \cline{1-1} \end{array}  .\]
This gives
 \[ M^{(3,1)}\cong \S^{(4)}\oplus \S^{(3,1)} .\]
For $\mu=(2,2)$, we have
 \[ \begin{array}{|c|c|c|c|} \hline 0&0&1&1 \\ \hline \end{array}  \quad\quad
 \begin{array}{|c|c|c|} \hline 0&0&1\\ \hline
                    1 \\ \cline{1-1} \end{array} \quad\quad
 \begin{array}{|c|c|} \hline 0&0\\ \hline
                   1& 1\\ \hline \end{array}  .\]
Thus we obtain
 \[ M^{(2,2)}\cong \S^{(4)}\oplus \S^{(3,1)} \oplus \S^{(2,2)} .\]
Finally, for $\mu=(2,1,1)$ we find
 \[ \begin{array}{|c|c|c|c|} \hline 0&0&1&2 \\ \hline \end{array}
 \quad\quad
 \begin{array}{|c|c|c|} \hline 0&0&1\\ \hline
                    2\\ \cline{1-1} \end{array} \quad\quad
  \begin{array}{|c|c|c|} \hline 0&0&2\\ \hline
                    1\\ \cline{1-1} \end{array} \quad\quad
  \begin{array}{|c|c|} \hline 0&0\\ \hline
                   1& 2 \\ \hline \end{array} \quad\quad
 \begin{array}{|c|c|} \hline 0&0 \\ \hline
                   1  \\ \cline{1-1}
                   2 \\ \cline{1-1} \end{array}
                   .\]
This means
 \[ M^{(2,1,1)}\cong \S^{(4)}\oplus  (\S^{(3,1)}\otimes \one_2) \oplus \S^{(2,2)} \oplus \S^{(2,1,1)}.\]
Putting this all together, we find the module decomposition of
$\pi(S_4)$ is given by

\begin{samepage} \begin{align}
 \pi(S_4)&= (M^{(4)}\otimes \one_3)\oplus (M^{(3,1)}\otimes
 \one_6)\oplus (M^{(2,2)}\otimes \one_3) \oplus
 (M^{(2,1,1)}\otimes \one_3) \nonumber \\
 &= (\S^{(4)}\otimes \one_3)\oplus \big((\S^{(4)}\oplus
 \S^{(3,1)})\otimes
 \one_6 \big) \label{dd-1} \\&\oplus \big( (\S^{(4)} \oplus \S^{(3,1)} \oplus \S^{(2,2)})
 \otimes \one_3 \big) \label{dd0} \\
 & \oplus
 \big((\S^{(4)}\oplus  (\S^{(3,1)}\otimes \one_2) \oplus \S^{(2,2)} \oplus \S^{(2,1,1)})\otimes \one_3\big) \label{dd}   \\
 &= (\S^{(4)}\otimes \one_{15})\oplus (\S^{(3,1)}\otimes
 \one_{15})
 \oplus (\S^{(2,2)}\otimes \one_6) \oplus (\S^{(2,1,1)}\otimes
 \one_3)
 . \label{dd1}
 \end{align}\end{samepage}
The direct sums in \eqref{dd-1}, \eqref{dd0} and \eqref{dd} are
understood to be `linked', as reflected in \eqref{dd1}.
It now follows that
 \begin{align}
 \fix(\E_\S) &= \pi(S_4)'' \nonumber \\ &\cong (\bbC\otimes\one_{15}) \oplus (\Ma_3\otimes \one_{15}) \oplus
 (\Ma_2\otimes \one_6) \oplus (\Ma_3 \otimes \one_3). \label{d3n4}
 \end{align}
Notice also that $\Ma_3$ is the largest full matrix algebra which
can be injected into $\fix(\E_\S)$ as a subalgebra.

Let us now describe the bases for the decomposition
\begin{eqnarray*}
\H_{3^4} &=& (\H_{4,0,0} \oplus \H_{0,4,0} \oplus\H_{0,0,4})
\oplus (\H_{2,1,1} \oplus\H_{1,2,1} \oplus\H_{1,1,2}) \\
& & \oplus (\H_{3,1,0} \oplus \H_{3,0,1} \oplus\H_{0,3,1} \oplus
\H_{1,3,0} \oplus\H_{0,1,3} \oplus\H_{1,0,3})
\end{eqnarray*}
which yields this algebra decomposition. This is easy for
$\la=(4)$. Indeed, for every $\H_{k_0,k_1,k_2}$ this is given by
the invariant average vector
 \[ h^{(4)}= \sum_{k_0(i_1,i_2,i_3,i_4)=k_0,\ldots,k_2(i_1,i_2,i_3,i_4)=k_2}  \ket{i_1i_2i_3i_4} . \]
In the following we will only discuss the case where $k_0\geq
k_1\geq k_2$ ({\it i.e.}, $\H_{4,0,0}$, $\H_{3,1,0}$,
$\H_{2,1,1}$). For $\la=\mu$, we have a natural embedding
$\S^{\la}\subseteq M^{\la}\cong \H_{k_0,k_1,k_2}$ given by
  \begin{equation}\label{form}
   h_t= \sum_{\sigma \in C_t} \sgn(\sigma)
   \sum_{(i_1,\ldots,i_4)\sim_t
  (i_1^t,\ldots,i_4^t)} \ket{i_1i_2i_3i_4},
  \end{equation}
for all $\la$-tableau $t$ of type $\mu$. Let us illustrate this in
our examples. If $\la=(3,1)$, we have $3$ standard tableaux
 \[ t_0=\begin{array}{|c|c|c|} \hline 1&2&3\\ \hline 4 \\ \cline{1-1} \end{array}  \quad\quad
 t_1=\begin{array}{|c|c|c|} \hline 1&2&4\\ \hline 3 \\ \cline{1-1} \end{array}  \quad\quad
  t_2=\begin{array}{|c|c|c|}  \hline 1&3&4\\ \hline 2 \\ \cline{1-1} \end{array}  .\]
The column stabilizer of $t_0$, $t_1$, $t_2$ is
$C_{t_0}=\{1,(14)\}$, $C_{t_1} = \{1, (13)\}$,
$C_{t_2}=\{1,(12)\}$. The space $\H_{(3,1)}$ has the basis
 \[ \ket{0001}, \quad \ket{0010},\quad \ket{0100},\quad \ket{1000}  .\]
Now, we define on $\H_{(3,1)}$
 \[ A_{t_i}= \sum_{\sigma\in C_{t_i}} \sgn(\sigma) \pi(\sigma) . \]
The range of $A_{t_i}$ is given by the vectors
 \begin{eqnarray*}
  h_{t_0} &=& \ket{0001}-\ket{1000},\\
  h_{t_1} &=& \ket{0010}-\ket{1000} ,\\
  h_{t_2} &=& \ket{0100}-\ket{1000} .
 \end{eqnarray*}
This provides us with the basis for
 \[ M^{(3,1)} \cong \H_{(3,1)}= \spn \{h^{(4)}\} \oplus \spn\{ h_{t_0},h_{t_1},h_{t_2}\}
  .\]

Now, we consider $\H_{(2,2,0)}$ spanned by
 \[ \ket{0011}, \quad \ket{0110}, \quad \ket{0101},\quad \ket{1100},\quad \ket{1010},\quad \ket{1001}
 .\]
For $\la=(3,1)$ we have the following list of $\la$-tableaux of
type $(2,2)$
  \[ \begin{array}{|c|c|c|} \hline 0&0&1\\ \hline 1 \\ \cline{1-1} \end{array}
  \,\,\,\,\,\,
 \begin{array}{|c|c|c|} \hline 0&1&1\\ \hline 0 \\ \cline{1-1} \end{array}
 \,\,\,\,\,\,
  \begin{array}{|c|c|c|} \hline 0&1&0\\ \hline 1 \\ \cline{1-1} \end{array}
  \,\,\,\,\,\,
  \begin{array}{|c|c|c|} \hline 1&1&0\\ \hline 0 \\ \cline{1-1} \end{array}
  \,\,\,\,\,\,
 \begin{array}{|c|c|c|} \hline 1&0&1\\ \hline 0 \\ \cline{1-1} \end{array} \,\,\,
 \,\,\,
  \begin{array}{|c|c|c|} \hline 1&0&0\\ \hline 1 \\ \cline{1-1} \end{array}
    .\]
Here we used $t_0={\Small \begin{array}{|c|c|c|} \hline 1&2&3\\
\hline 4
\\ \cline{1-1}
\end{array}}$. Only the first tableaux is semistandard and yields
an injection
 \[ \Theta\equiv \Theta_{{\Small \begin{array}{|c|c|c|} \hline 0&0&1\\ \hline 1 \\ \cline{1-1} \end{array}}}: \S^{(3,1)}\rightarrow \H_{(2,2)}
\]
with
\[
 \Theta(e_{\{t_0\}})=
 \sum_{(i_1,\ldots,i_4)\sim_{t_0} (0,0,1,1)} \ket{i_1i_2i_3i_4}. \]
This means
 \[ \Theta(e_{\{t_0\}})= \ket{0011}+ \ket{0101}+ \ket{1001}. \]
(See the example in the last section for $n=5$.) Further, this
vector is a cyclic vector for the image of $\S^{(3,1)}$ in
$\H_{(2,2)}$. The polytabloid is
 \[ e_{t_0}= \sum_{\sigma \in C_{t_0}} \sgn(\sigma) \sigma e_{\{t_0\}}
  = e_{\{t_0\}}-(14)e_{\{t_0\}} .\]
Therefore $\Theta(\S^{(3,1)})$ is the module generated by
 \[ \Theta(e_{t_0})= h_{t_0} = \ket{0011}- \ket{1010}+
 \ket{0101}- \ket{1100}  .\]
Equivalently,
 \begin{eqnarray*}
  \Theta(\S^{(3,1)})&=& \spn\{ h_{t_0}, (12)h_{t_0}, (13) h_{t_0}
  \} \\
  &=&  \spn\big\{\ket{0011}-\ket{1010}+ \ket{0101}-\ket{1100}, \\
 & & \ket{1100}-\ket{1001}+ \ket{0110}-\ket{0011},  \\
 & &\ket{0110}-\ket{0101}+ \ket{1010}-\ket{1001}\big\}  .
 \end{eqnarray*}
Another way to find a basis is to consider $t_1={\Small
\begin{array}{|c|c|c|} \hline 1&2&4\\ \hline 3 \\ \cline{1-1}
\end{array}}$. In this case,
$C_{t_1}=\{1,(13)\}$,
 \[
 \Theta(e_{\{t_1\}})
 = \ket{0011}+ \ket{0110}+\ket{1010}
  \]
and
 \[ h_{t_1}= \Theta(e_{\{t_1\}})-(13)\Theta(e_{\{t_1\}})=
 \ket{0011}-\ket{1001}+ \ket{0110}- \ket{1100}  .\]
Similarly for $t_2={\Small \begin{array}{|c|c|c|} \hline 1&3&4\\
\hline 2 \\ \cline{1-1}
\end{array}}$, we have $C_{t_1} =\{1,(12)\}$ and
 \[ \Theta(e_{\{t_2\}}) =
  \ket{0101}+\ket{0110}+ \ket{1100}
  \]
and
 \[ h_{t_2}= \Theta(e_{\{t_2\}})-(12)\Theta(e_{\{t_2\}})=
 \ket{0101}-\ket{1001}+ \ket{0110}- \ket{1010} .\]

The copy of $\S^{(2,2)}$ in $\H_{(2,2)}$ is again easy to find. We
recall that $\S^{(2,2)}$ is spanned by the standard tableaux
$\{e_{s_0}, e_{s_1}\}$ where
 \[ s_0=\begin{array}{|c|c|} \hline 1&2\\ \hline 3&4 \\ \hline \end{array}  \quad \quad
 s_1=\begin{array}{|c|c|} \hline 1&3\\ \hline 2&4 \\ \hline \end{array}  .\]
The column stabilizers are given by
$C_{s_0}=\{1,(13),(24),(13)(24)\}$ and by $C_{s_1}
=\{1,(12),(34),(12)(34)\}$. This yields operators on $\H_{(2,2)}$,
\[ A_{s_0}= 1-\pi((13))-\pi((24))+ \pi((13)(24))  \]
and
\[ A_{s_1}= 1-\pi((12))-\pi((34))+ \pi((12)(34)) . \]
Applied to the unit vectors, we find the ranges
 \[ \ran (A_{s_0}) = \ket{0110}-\ket{1010}-\ket{0101}+\ket{1001}\]
and
 \[ \ran(A_{s_1}) = \ket{0011}-\ket{1001}-\ket{0110}+\ket{1100}. \]

Finally, we consider $\H_{(2,1,1)}$ with basis
 \begin{eqnarray*}
 \big\{\ket{0012},\ket{0021},\ket{0102},\ket{0120},\ket{0201},\ket{0210},\ket{1002},
 \ket{1020},  \ket{1200}, \\ \ket{2001},\ket{2010},\ket{2100}\big\}
  .\end{eqnarray*}
The representation of $\S^{(4)}$ is $1$-dimensional, given by the
average of all these vectors. There are two copies of $\S^{(3,1)}$
corresponding to the two $(3,1)$-tableaux of type $(2,1,1)$
 \[ T= \begin{array}{|c|c|c|} \hline 0&0&1\\ \hline 2 \\ \cline{1-1} \end{array}  \quad\quad
 T'= \begin{array}{|c|c|c|} \hline 0&0&2\\ \hline 1 \\ \cline{1-1} \end{array}  .\]
The basis for $\S^{(3,1)}$ is given by $\{ e_{t_0}, e_{t_1},
e_{t_2}\}$ where
 \[ t_0= \begin{array}{|c|c|c|} \hline 1&2&3\\ \hline 4 \\ \cline{1-1} \end{array}  \quad
 t_1 = \begin{array}{|c|c|c|} \hline 1&2&4\\ \hline 3 \\ \cline{1-1} \end{array}  \quad
 t_2= \begin{array}{|c|c|c|} \hline 1&3&4\\ \hline 2 \\ \cline{1-1} \end{array}  .\]
Following the definition of $\Theta_{T}(e_{t_j})$, we get
 \[  \gamma_{t_0}(T) = (0,0,1,2), \quad \gamma_{t_1}(T)= (0,0,2,1)
 ,\quad
 \gamma_{t_2}(T)= (0,2,0,1)  .\]
Using  row equivalence, we are allowed to permute the entries
$\{1,2,3\}$ for $t_0$, the entries $\{1,2,4\}$ for $t_1$ and
$\{1,3,4\}$ for $t_2$ and thus
 \begin{align*}
 \Theta_{T}(e_{\{t_0\}})
 &=\ket{0012}+\ket{0102}+\ket{1002}, \\
 \Theta_{T}(e_{\{t_1\}})
  &= \ket{0021}+\ket{0120}+\ket{1020},  \\
 \Theta_{T}(e_{\{t_2\}})
  &= \ket{0201}+\ket{0210}+\ket{1200} .
  \end{align*}
For $t_0$, $t_1$, $t_2$ we have to apply, respectively,
$A_{T,t_0}=\one -\pi(14)$, $A_{T,t_1}=\one-\pi(13)$ and
$A_{T,t_2}=\one-\pi(12)$ in order to obtain the image of the
polytabloids:
 \begin{align*}
 h_{T,t_0} &=
 \ket{0012}-\ket{2010}+\ket{0102}-\ket{2100}+\ket{1002}-\ket{2001} ,\\
 h_{T,t_1} &=
 \ket{0021}-\ket{2001}+\ket{0120}-\ket{2100}+\ket{1020}-\ket{2010} ,\\
 h_{T,t_2} &=
 \ket{0201}-\ket{2001}+\ket{0210}-\ket{2010}+\ket{1200}-\ket{2100} .
 \end{align*}
This is our first copy of $\S^{(3,1)}$. For the second, we
exercise the same procedure in the case of $T'$.
 \begin{align*}
 \Theta_{T'}(e_{\{t_0\}})
 &=\ket{0021}+\ket{0201}+\ket{2001},\\
 \Theta_{T'}(e_{\{t_1\}})
  &= \ket{0012}+\ket{0210}+\ket{2010}, \\
 \Theta_{T'}(e_{\{t_2\}})
  &= \ket{0102}+\ket{0120}+\ket{2100} .
  \end{align*}
This provides us with
 \begin{align*}
 h_{T',t_0} &=
 \ket{0021}-\ket{1020}+\ket{0201}-\ket{1200}+\ket{2001}-\ket{1002} ,\\
 h_{T',t_1} &=
 \ket{0012}-\ket{1002}+\ket{0210}-\ket{1200}+\ket{2010}-\ket{1020} ,\\
 h_{T',t_2} &=
 \ket{0102}-\ket{1002}+\ket{0120}-\ket{1020}+\ket{2100}-\ket{1200} .
 \end{align*}
We have one copy of $\S^{(2,2)}$ which is spanned by
 \[ s_0= \begin{array}{|c|c|} \hline 1&2\\ \hline 3&4 \\ \hline \end{array}  \quad\quad
  s_1= \begin{array}{|c|c|} \hline 1&3\\ \hline 2&4 \\ \hline \end{array}  .\]
Our $(2,2)$ tableau of type $(2,1,1)$ is given by $T={\Small
\begin{array}{|c|c|} \hline 0&0\\ \hline 1&2 \\ \hline \end{array}}$. This gives
 \[ \Theta(e_{\{s_0\}})= \ket{0012}+\ket{0021} \quad \mbox{and} \quad
  \Theta(e_{\{s_1\}})= \ket{0102}+\ket{0201} . \]
The operator is $\one-\pi(13)-\pi(24)+\pi((13)(24))$, determined
by $C_{s_0} = \{ 1, (13), (24), (13)(24)\}$, and thus
 \[ h_{s_0}= \ket{0012}-\ket{1002}-\ket{0210}+\ket{1200}+
 \ket{0021}-\ket{2001}- \ket{0120}+\ket{2100} \]
and similarly  for $s_1$ we apply $\one
-\pi(12)-\pi(34)+\pi((12)(34))$ to obtain
 \[ h_{s_1}= \ket{0102}-\ket{1002}-\ket{0120}+\ket{1020}+\ket{0201}-\ket{2001}-\ket{0210}+\ket{2010} . \]

Finally we consider the copy of $\S^{(2,1,1)}$, which has  basis
$\{e_{r_0}, e_{r_1}, e_{r_2}\}$ where
 \[ r_0= \begin{array}{|c|c|} \hline 1&2\\ \hline 3\\ \cline{1-1} 4 \\ \cline{1-1} \end{array} \quad\quad
 r_1= \begin{array}{|c|c|} \hline 1&3\\ \hline 2\\ \cline{1-1} 4 \\ \cline{1-1} \end{array} \quad\quad
 r_2= \begin{array}{|c|c|} \hline 1&4\\ \hline 2\\ \cline{1-1} 3 \\ \cline{1-1} \end{array}
 .  \]
Here $T={\Small \begin{array}{|c|c|} \hline 0&0\\
\hline 1\\ \cline{1-1} 2 \\ \cline{1-1} \end{array}}$. This yields
 \[ \Theta_T(e_{\{r_0\}})= \ket{0012} , \quad
 \Theta_T(e_{\{r_1\}}) = \ket{0102},\quad \Theta_T(e_{\{r_2\}})= \ket{0120}
 .\]
The column stabilizer of $r_2$ is given by all permutations which
leave $\{1,2,3\}$ invariant. This yields
 \begin{align*}
  h_{r_2} &=    \ket{0120}
  -\ket{0210}-\ket{1020}+\ket{1200}+\ket{2010}-\ket{2100} .
 \end{align*}
Similarly, we have to look for all permutations of $\{1,3,4\}$ in
the column stabilizer of $r_0$ and we obtain
 \begin{align*}
  h_{r_0} &=  \ket{0012}
  -\ket{0021}-\ket{1002}+\ket{1020}+\ket{2001}-\ket{2010}.
  \end{align*}
For the column stabilizer of $r_1$, we may permute $\{1,2,4\}$ and
hence
 \begin{align*}
  h_{r_1} &=  \ket{0102}
  -\ket{0201}-\ket{1002}+\ket{1200}+\ket{2001}-\ket{2100} .
  \end{align*}

By equation (\ref{d3n4}), the largest full matrix algebra $\Ma_k$
that can be injected into the noise commutant here is $\Ma_3$,
identified with the subalgebras of $\A_\S^\prime$ unitarily
equivalent to either $ \one_{15} \otimes \Ma_3$ or $ \one_3
\otimes \Ma_3$. Let us explicitly identify the copy of $\one_3
\otimes \Ma_3 $. The set $\{ h_{r_0}, h_{r_1}, h_{r_2} \}$ yields
the copy of $S^{(2,1,1)}$ inside $\H^{(2,1,1)}$. A similar
analysis yields the basis for the copy of $S^{(2,1,1)}$ inside
$\H^{(1,2,1)}$. It is generated by $T'={\Small \begin{array}{|c|c|} \hline 0&1\\
\hline 1\\ \cline{1-1} 2 \\ \cline{1-1} \end{array}}$ and in this
case
 \begin{align*}
 \Theta_{T'}(e_{\{r_0\}})
 &=\ket{0112}+\ket{1012}, \\
 \Theta_{T'}(e_{\{r_1\}})
  &= \ket{0112}+\ket{1102},  \\
 \Theta_{T'}(e_{\{r_2\}})
  &= \ket{0121}+\ket{1120} .
  \end{align*}
Thus we have
\[
\left\{ \begin{array}{rcl}
  h_{r_0}' &=&  \ket{0112}
  -\ket{0121}-\ket{1102}+\ket{1120}+\ket{2101}-\ket{2110} \\
h_{r_1}' &=&  \ket{0112}
  -\ket{0211}-\ket{1012}+\ket{1210}+\ket{2011}-\ket{2110} \\
 h_{r_2}' &=&    \ket{0121}
  -\ket{0211}-\ket{1021}+\ket{1201}+\ket{2011}-\ket{2101}
  \end{array}\right..
\]
Further, the basis for the copy of $S^{(2,1,1)}$ inside
$\H^{(1,1,2)}$ is generated by $T''={\Small \begin{array}{|c|c|} \hline 0&2\\
\hline 1\\ \cline{1-1} 2 \\ \cline{1-1} \end{array}}$ and in this
case
 \begin{align*}
 \Theta_{T''}(e_{\{r_0\}})
 &=\ket{0212}+\ket{2012}, \\
 \Theta_{T''}(e_{\{r_1\}})
  &= \ket{0122}+\ket{2102},  \\
 \Theta_{T''}(e_{\{r_2\}})
  &= \ket{0122}+\ket{2120} .
  \end{align*}
Thus we have
\[
\left\{ \begin{array}{rcl}
  h_{r_0}'' &=&  \ket{0212}
  -\ket{0221}-\ket{1202}+\ket{1220}+\ket{2201}-\ket{2210} \\
h_{r_1}'' &=&  \ket{0122}
  -\ket{0221}-\ket{1022}+\ket{1220}+\ket{2021}-\ket{2120} \\
 h_{r_2}'' &=&    \ket{0122}
  -\ket{0212}-\ket{1022}+\ket{1202}+\ket{2012}-\ket{2102}
  \end{array}\right..
\]

Let $P_{(2,1,1)}$ be the projection of $\H$ onto the span of $\{
h_{r_i}, h_{r_j}', h_{r_k}'': 0 \leq i,j,k \leq 2\}$. Then
$P_{(2,1,1)}$ is a minimal central projection for $\A_\S^\prime$
and the `compression subalgebra' $P_{(2,1,1)} \A_\S^\prime
P_{(2,1,1)} = \A_\S^\prime P_{(2,1,1)} \subset \A_\S^\prime$ is
unitarily equivalent to $\one_3 \otimes \Ma_3$. In fact, with
respect to the ordered basis
\[
\{ h_{r_0}, h_{r_1}, h_{r_2}, h_{r_0}', h_{r_1}', h_{r_2}',
h_{r_0}'', h_{r_1}'', h_{r_2}'' \}
\]
for $P_{(2,1,1)}\H$, we have the matrix representations
\[
\A_\S^\prime P_{(2,1,1)} = \left\{ \left( \begin{matrix} A & 0 & 0 \\
 0 & A & 0 \\ 0 & 0 & A \end{matrix}\right) : A\in\Ma_3 \right\}.
 \]
Note that the subspaces spanned by the sets $\{h_{r_i}\}$,
$\{h_{r_i}'\}$ and $\{h_{r_i}''\}$ are perpendicular, but the
vectors within each of these sets do not form an orthogonal basis
for the corresponding subspace.


\section{Conclusion}\label{S:conclusion}

We have investigated the operator algebras of fixed points for a
class of quantum channels we call universal collective rotation
channels $\E_\S$. This class includes as a subclass the well-known
class of collective rotation channels. We showed that such
channels always have an abundance of noiseless subsystems and gave
a method for explicitly computing them. In particular, the Young
tableaux machine gives a clean approach for this process. In lower
dimensional cases ({\it e.g.} when $d=2$), our approach is more
technical when compared to others in the literature (for instance
\cite{HKLP}). However, an important advantage of the Young
tableaux approach for higher dimensional cases is that it is
particularly amenable to computations.

An issue we have not pursued here concerns the channels generated
by non-maximal sets $\S$. The $d^n$-dimensional representations of
$J_x,J_y,J_z$ considered in \cite{HKLP} provide such an example,
but we would expect there to be other interesting non-trivial
examples of channels $\E_\S$ for non-maximal $\S$. We emphasize
that even for non-maximal $\S$ there is an abundance of noiseless
subsystems because  $\fix(\E_\S)$ contains $\pi(S_n)''$.

We also wonder what other representations of $S_n$ correspond to
physically meaningful unital  channels, beyond $\pi$ and its
subrepresentations (which correspond to the compressions of
ucr-channels). The recent preprint \cite{BCH} of Bacon, et al,
appears to present further insights into this topic, and also
shows how the unitary base change from the standard basis to the
basis given by the Young tableaux can be efficiently computed
using quantum circuits.

\vspace{0.1in}


{\noindent}{\it Acknowledgements.} We are grateful to the referee
for helpful comments. We would like to thank John Holbrook,
Raymond Laflamme and David Poulin for enlightening conversations.
The first author was partially supported by NSF grant DMS
03-01116. The second and third authors were partially supported by
NSERC grants. The third author also gratefully acknowledges
support from the Institute for Quantum Computing and the Perimeter
Institute for Theoretical Physics.


\vspace{0.1in}




\begin{thebibliography}{99}


\bibitem{AB} D. Aharonov,  M. Ben-Or, {\it Fault-tolerant quantum
computation with constant error,} In Proc. 29th. Ann. ACM Symp. on
Theory of Computing, page 176, New York, 1998, ACM.
arxiv.org/quant-ph/9906129, quant-ph/9611025.









\bibitem{Arvinvite} W. Arveson,
\textit{An invitation to $\ca$-algebras}, Graduate Texts in
Mathematics, No. 39, Springer-Verlag, New York-Heidelberg, 1976.


\bibitem{BCH} D. Bacon, I.L. Chuang, A.W. Harrow,
\textit{Efficient quantum circuits for Schur and Clebsch-Gordon
transforms}, arxiv.org/quant-ph/0407082.


\bibitem{BRS2} S.D. Bartlett, T. Rudolph, R.W. Spekkens,
\textit{Decoherence-full subsystems and the cryptographic power of
a private shared reference frame}, arXiv.org/quant-ph/0403161.

\bibitem{BRS} S.D. Bartlett, T. Rudolph, R.W. Spekkens,
\textit{Classical and quantum communication without a shared
reference frame}, Phys. Rev. Lett. {\bf 91}, 027901 (2003).








\bibitem{LPS} J.-C. Boileau, D. Gottesman, R. Laflamme, D. Poulin,
R.W. Spekkens, \textit{Robust polarization-based quantum key
distribution over collective-noise channel},
arXiv.org/quant-ph/0306199.











\bibitem{BS} P. Busch, J. Singh,
\textit{Luders theorem for unsharp quantum effects,} Phys. Lett.
A, \textbf{249} (1998), 10-24.


\bibitem{Choi} M.D. Choi,
\textit{Completely positive linear maps on complex matrices,} Lin.
\ Alg.\ Appl. \ \textbf{10} (1975), 285-290.


\bibitem{QM12} C. Cohen-Tannoudji, B. Diu, F. Laloe,
\textit{Quantum Mechanics, Volume One \& Two}, John Wiley \& Sons,
Toronto, 1977.







\bibitem{byeg} K.R. Davidson,
\textit{$\ca$-algebras by example}, Fields Institute Monographs,
{\bf 6}, Amer. Math. Soc., Providence, 1996.














\bibitem{DG} L.-M. Duan, G.-C. Guo,
\textit{Preserving coherence in quantum computation by pairing
quantum bits}, Phys. Rev. Lett. {\bf 79} (1997), 1953.









\bibitem{Fi} S. De Filippo,
\textit{Quantum computation using decoherence-free states of the
physical operator algebra,} Phys. Rev. A {\bf 62}, 052307 (2000).


\bibitem{Lnoiseless} E.M. Fortunato, L. Viola, M.A. Pravia, E. Knill, R.
Laflamme, T.F. Havel, D.G. Cory, \textit{Exploring noiseless
subsystems via nuclear magnetic resonance,} Phys. Rev. A {\bf 67},
062303 (2003).




\bibitem{FultonH} W. Fulton, J. Harris,
\textit{Representation theory, a first course}, Springer-Verlag
New York, 1991.







\bibitem{Green} J.A. Green,
\textit{Polynomial representations of $GL_n$}, Lecture Notes in
Mathematics 830, Springer-Verlag, New York, 1980.





\bibitem{HKLP} J.A. Holbrook, D.W. Kribs, R. Laflamme, D. Poulin,
\textit{Noiseless subsystems for collective rotation channels in
quantum information theory}, Integral Equations \& Operator
Theory, to appear.


\bibitem{HKL} J.A. Holbrook, D.W. Kribs, R. Laflamme,
\textit{Noiseless subsystems and the structure of the commutant in
quantum error correction}, Quantum Information Processing {\bf 2}
(2003), 381-419.










\bibitem{James} G. James, A. Kerber,
\textit{The representation theory of the symmetric group,}
Encyclopedia of Mathematics and its Applications, Addison-Wesley
Publishing Company, Toronto, 1981.




\bibitem{Junge} M. Junge,
\textit{Noncommutative Poisson random measure,} preprint, 2003.

\bibitem{KBLW} J. Kempe, D. Bacon, D.A. Lidar, K.B. Whaley,
\textit{Theory of decoherence-free fault-tolerant universal
quantum computation,} Phys. Rev. A {\bf 63}, 042307 (2001).




\bibitem{IntroQEC} E. Knill, R. Laflamme, A. Ashikhmin, H. Barnum, L.
Viola, W.H. Zurek, \textit{Introduction to Quantum Error
Correction}, Los Alamos Science, November 27, 2002.

\bibitem{KLV} E. Knill, R. Laflamme, L. Viola,
\textit{Theory of quantum error correction for general noise},
Phys. Rev. Lett. {\bf 84} (2000), 2525-2528.

\bibitem{KLZ} E. Knill, R. Laflamme, W. H. Zurek,
\textit{Resilient quantum computation: error models and
thresholds,} Science {\bf 279} (1998), 342-345.

\bibitem{KL97} E. Knill, R. Laflamme,
\textit{A theory of quantum error-correcting codes,} Phys. Rev. A
{\bf 55} (1997), 900.

\bibitem{Kraustext} K. Kraus,
\textit{States, Effects and Operations: Fundamental Notions of
Quantum Theory,} Lecture Notes in Physics, vol. 190, Berlin:
Springer-Verlag, 1983.

\bibitem{Kraus} K. Kraus,
\textit{General state changes in quantum theory,} Ann. \ Physics \
\textbf{64} (1971), 311-335.



\bibitem{Kquantumsurvey} D.W. Kribs,
\textit{A quantum information theory primer for operator
theorists}, preprint, 2004.

\bibitem{Kchannel} D.W. Kribs,
\textit{Quantum channels, wavelets, dilations, and representations
of $\O_n$}, Proc. Edin. Math. Soc., {\bf 46} (2003).









\bibitem{Rayprivate} R. Laflamme, private communication.


\bibitem{LCW} D. A. Lidar, I. L. Chuang, K. B. Whaley,
\textit{Decoherence free subspaces for quantum computation}, Phys.
Rev. Lett.  {\bf 81} (1998), 2594.








\bibitem{NC} M.A. Nielsen, I.L. Chuang,
\textit{Quantum computation and quantum information}, Cambridge
University Press, 2000.



\bibitem{Paulsentext} V. Paulsen,
\textit{Completely bounded maps and operator algebras,} Cambridge
University Press, Cambridge, United Kingdom, 2002.















\bibitem{Preskillplenum} J. Preskill,
\textit{Reliable quantum computers,} Proc. R. Soc. Lond. A, {\bf
454} (1998), 385-410.













\bibitem{Stine} W.F. Stinespring,
\textit{Positive functions on $\ca$-algebras}, Proc. Amer. Math.
Soc. {\bf 6} (1955), 211-216.

\bibitem{Tak} M. Takesaki,
\textit{Theory of operator algebras I}, Springer-Verlag, New
York-Heidelberg, 1979.

\bibitem{VKL} L. Viola, E. Knill, R. Laflamme,
\textit{Constructing qubits in physical systems}, J. Phys. A {\bf
34}, 7067 (2001).

\bibitem{VFPKLC} L. Viola, E.M. Fortunato,  M.A. Pravia, E. Knill, R.
Laflamme, D.G. Cory, \textit{Experimental realization of noiseless
subsystems for quantum information processing,} Science {\bf 293},
2059 (2001).




\bibitem{ZL} P. Zanardi, S. Lloyd,
\textit{Topological protection and quantum noiseless subsystems},
Phys. Rev. Lett. {\bf 90}, 067902 (2003).

\bibitem{Zan} P. Zanardi,
\textit{Stabilizing quantum information}, Phys. Rev. A {\bf 63},
012301 (2001).

\bibitem{ZR} P. Zanardi,  M. Rasetti,
\textit{Noiseless quantum codes}, Phys. Rev. Lett. {\bf 79}
(1997), 3306.

\bibitem{Zar} V. Zarikian,
\textit{Algorithms for operator algebra calculations}, Lin. Alg.
Appl., to appear.


\end{thebibliography}
\end{document}